\documentclass[3p,times]{elsarticle}

\usepackage{enumitem}
\usepackage{hyperref}
\usepackage{amsmath}
\usepackage{pifont}
\usepackage[utf8]{inputenc}
\makeatletter
\def\LaTeX{\leavevmode L\raise.42ex
    \hbox{\kern-.3em\size{\sf@size}{0pt}\selectfont A}\kern-.15em\TeX}
\makeatother

\newcommand{\BibTeX}{{\rm B\kern-.05em{\sc
          i\kern-.025emb}\kern-.08em\TeX}}

\makeatletter
\def\@currentlabel{2.1}\label{e:dispaa}
\def\@currentlabel{2.21}\label{e:dispau}
\def\@currentlabel{2.22}\label{e:dispav}
\def\@currentlabel{2.23}\label{e:dispaw}
\def\@currentlabel{2.24}\label{e:dispax}
\def\theequation{\thesection.\@arabic\c@equation}
\makeatother

\renewcommand{\theequation}{\arabic{section}.\arabic{equation}}

\newcommand{\R}{\mathbb R}

\newcommand{\N}{\mathbb N}

\def \O{\Omega}
\everymath{\displaystyle}
\newtheorem{theorem}{Theorem}

\newtheorem{thm}{Theorem} [section]
\newtheorem{lem}{Lemma} [section]
\newtheorem{prop}{Proposition} [section]

\newtheorem{rem}{Remark}[section]

\renewcommand{\theequation}{\thesection.\arabic{equation}}
\renewcommand{\thesection}{\arabic{section}}
\renewcommand{\theequation}{\thesection.\arabic{equation}}
\let\ssection=\section\renewcommand{\section}{\setcounter{equation}{0}\ssection}
\begin{document}
\begin{frontmatter}

\title{High-order Kirchhoff problems in bounded and unbounded domains}
\author[mk1,mk2]{Hamdani Mohamed Karim}
\ead{hamdanikarim42@gmail.com}
\cortext[cor1]{Corresponding author:Tel:+21698688775.}
\author[ah3,ah1,ah2]{Harrabi Abdellaziz\corref{cor1}}
\ead{abdellaziz.harrabi@yahoo.fr}
\begin{center}
\address[mk1]{ Mathematics Department, University of Sfax, Faculty of Science of Sfax, Sfax, Tunisia.}
\address[mk2] {Military School of Aeronautical Specialities, Sfax, Tunisia.}
\address[ah3]{Mathematics Department, University of Kairouan, Higher Institute of Applied Mathematics and Informatics, Kairouan, Tunisia.}
\address[ah1]{Mathematics Department, Northern Borders University, Arar, Saudi Arabia.}
\address[ah2]{Senior associate in the Abdus Salam International Centre for Theoretical Physics, Trieste, Italy.}
\end{center}
\begin{abstract}
  Consider the following $m-$polyharmonic Kirchhoff problem:
\begin{eqnarray} \label{ea}
\begin{cases}
M\left(\int_{\O}|D_r u|^{m} +a|u|^m\right)[\Delta^r_m u +a|u|^{m-2}u]= K(x)f(u) &\mbox{in}\quad \Omega, \\
u=\left(\frac{\partial}{\partial \nu}\right)^k u=0, \quad &\mbox{on}\quad \partial\Omega, \quad k=1, 2,.....  , r-1,
\end{cases}
\end{eqnarray}
where $r \in \N^*$, $m >1$, $N\geq rm+1$, $a\geq 0$, $K\in L^{\infty}(\O)$ is a positive weight function, $M \in C([0,+\infty))$ and $f\in
C(\mathbb{R})$ which will be specified later. We will study problem \eqref{ea} in the following different type of domains:
\begin{enumerate}
  \item $a=0$ and $K\in L^{\infty}(\O)$ is a positive weight function if $\Omega$ is a smooth bounded domain of $\R^N$.
  \item $a>0$ and $K\in L^{\infty}(\O)\cap L^{p}(\O)$, $p \geq 1$ if $\Omega$  is an unbounded smooth domain.
  \item $\O=\R^N$ and $a=0$ (which called the $m\gamma$-zero mass case).
\end{enumerate}

 We prove the  existence of infinitely many solutions of \eqref{ea} for some odd functions $f$ in $u$ satisfying subcritical growth conditions at infinity which are weaker than the analogue of the Ambrosetti-Rabinowitz condition and the standard subcritical polynomial growth. The new aspect consists in employing the Schauder basis of $W_0^{r,m}(\O)$ to verify the geometry of the symmetric mountain pass theorem without any control on $f$ near $0$ if $\Omega$ is a bounded domain and under a suitable condition at $0$ if $\Omega$ is a unbounded domain allowing only to derive the variational setting of \eqref{ea}. Moreover, we introduce a positive quantity $\lambda_M$ similar to the first eigenvalue of the $m$-polyharmonic operator to find a mountain pass solution.
\end{abstract}
\begin{keyword}
Palais-Smale condition\sep Symmetric mountain pass theorem\sep Schauder basis \sep Krasnoselskii genus theory \sep $m$-polyharmonic operator\sep  Kirchhoff equations\sep Zero mass case.
\PACS{Primary: 35J55, 35J65; Secondary: 35B65.}
\end{keyword}
\end{frontmatter}

\section{Introduction}\label{nnn}
This paper is concerned with the existence and multiplicity of solutions to the following $m-$polyharmonic Kirchhoff problem:
\begin{eqnarray} \label{10}
\begin{cases}
M\left(\int_{\O}|D_r u|^{m} +a|u|^m\right)[\Delta^r_m u +a|u|^{m-2}u]= K(x)f(u) &\mbox{in}\quad \Omega, \\
u=\left(\frac{\partial}{\partial \nu}\right)^k u=0, \quad &\mbox{on}\quad \partial\Omega, \quad k=1, 2,.....  , r-1,
\end{cases}
\end{eqnarray}
where $r \in \N^*$, $m >1$, $N\geq rm+1$, $a\geq 0$, $K\in L^{\infty}(\O)$ is a positive weight function, $M \in C([0,+\infty))$ and $f\in
C(\mathbb{R})$ which will be specified later. The $m-$polyharmonic operator $\Delta^r_m$ is defined by

$$\Delta^r_m u=\begin{cases}-div \left\{\Delta^{j-1}(|\nabla \Delta^{j-1} u|^{m-2} \nabla \Delta^{j-1} u)\right\}, \text{ if $r=2j-1$}\\
\Delta^j (|\Delta^{j} u|^{m-2} \Delta^{j} u), \text{ if $r=2j$}
\end{cases}j\in \N^*,$$
which becomes the usual polyharmonic operator for $m=2$, namely $(-\Delta )^r$. Define the main $r-$order differential operator by
$$\textit{D}_r u=\begin{cases}
 \nabla\Delta^{j-1}u & \text{if $r=2j-1$},\\
\Delta^j u & \text{if $r=2j$},
\end{cases}\;\;j\in \N^*.
 $$
Note that $D_r u$ is an $N-$vectorial operator when $r$  is odd and $N > 1$, while it is a scalar operator when $r$ is even. Denote $E=W_0^{r,m}(\O)$   endowed with the following norm:
$$\|u\|^m=\int_{\R^N}(|D_r u|^{m} +a|u|^m) \; \mbox{ where } a=0 \mbox{ if } \O \mbox{ is bounded}.$$
The adequate functional space corresponding to the $m\gamma$-zero mass case is $E= D^{r,m}(\mathbb{ R^{N}})$, defined as the completion of $C_c^\infty(\R^N)$ with respect the norm $\|u\|^m=\int_{\R^N}|D_r u|^{m}$.

 $(E,\|.\|)$ is a separable, uniformly convex, reflexive, real Banach space. Denote $E^*$ the  dual space of $E$, $p^* =\frac{mN}{N-rm}$ the Sobolev critical exponent and $q^*=\frac{p^*}{p^*-1}$ is the conjugate exponent of $p^*$. Recall also the Gagliardo-Nirenberg-Sobolev inequality which will be used a number of times later
\begin{eqnarray}\label{Gagliardo}
||u||_{L^{p^*}(\O) }\leq C\|u\|, \forall u \in E.
\end{eqnarray}

In recent years, there has been an increasing interest in studying problem \eqref{10}, which has a broad background in many different applications, such as game theory, mathematical finance, continuum mechanics,
phase transition phenomena, population dynamics and minimal surface. The reader may consult \cite{AM,BR,BRS,CP,CWL,HZ4,K,L,LLS1,P,RXZ1,RXZ2} and the references therein. Most part of these papers deal only with the non-degenerate case, that is when $M(\tau)\geq s>0$ for all $\tau\in [0,+\infty) $.

When $M\equiv1$, $m=2$, $r=1$, $a = 0$ and $\O=\R^\mathbb{N}$, then \eqref{10} boils down to the following elliptic
problem with zero mass:
\begin{eqnarray}\label{mm=1}
-\Delta u =K(x) f(u)\;\;
\mbox{ in }\R^N,
\end{eqnarray}
for which some  existence results have been established (see for example \cite{ASM,AZ}).

There are few papers considering Kirchhoff type problems on $\R^N$ (see \cite{FF1,JW,LLS0,LLS1}). Recently, Li-Li-Shi employed variational method to obtain the existence of positive solutions to the following zero mass problem

\begin{eqnarray}\label{prob-kirchh}
\left(a+\lambda\int_{\R^N}|\nabla u|^2+\lambda b \int_{\R^N}u^2\right)[-\Delta u+bu]= K(x)f(u), \; \mbox{ in }\; \R^N,
 \end{eqnarray}
where $b=0$, $\lambda\geq 0$ and $K$ is a nonnegative weight function satisfying
 $K \in [L^p(\R^N) \cap L^\infty(\R^N)]\setminus \{0\}$ for some $p\geq \frac{2N}{(N + 2)}.$ In particular, they used a cut-off functional and Pohozaev type identity to obtain the bounded Palais–Smale sequences \cite{LLS1}.  Similar result has been obtained in \cite{LLS0} related to the delicate case $K(x)\equiv 1$, $N\geq3$, $a, b$ are positive constants, $\lambda\geq0$.  The
sublinear case was also discussed by Feng-Feng in \cite{FF1} to  derive multiplicity results by using an extension of Clark's theorem.

The study of Kirchhoff type equations has already been extended to the case involving the $p$-Laplacian.  Very recently, Han-Ma-He \cite{HMH} considered the following problem with the zero mass case
\begin{eqnarray*}
\begin{cases}
-\left(a+b\int_{\R^N}|\nabla u|^pdx\right)\Delta_p u+V(x)|u|^{p-2}u= K(x)f(u), \; \mbox{ in }\; \R^N.

\end{cases}
 \end{eqnarray*}
 Using Nehari manifold method, they obtained the existence of least energy sign-changing solution (see also \cite{Hamdani}).
\subsection{\bf Main Assumptions and the variational setting of \eqref{10}}
We will study \eqref{10} related to the following different type of domains:
\begin{enumerate}
  \item $\Omega$ is a smooth bounded domain and $a=0$.
  \item $\Omega$  is an unbounded smooth domain and $a>0$.
  \item $\O=\R^N$ and $a=0$ (which called the $m\gamma$-zero mass case).
\end{enumerate}
In all cases, we assume that $K \in L^\infty(\O)$ is a positive weight function. However, if $\O$ is {\bf unbounded }, we need in addition:\\
${\mathbf{(K_p)}}$: there exists $\mathbf{p\geq 1}$ such that $K \in  L^{p}(\O)$.

We relax the global structural assumption assumed in \cite{CP} into the following:\\
$(M_1):$  $M: [0,+\infty) \rightarrow [0,+\infty)$ a continuous function satisfying: there exist $\tau_0\geq0$ and $\gamma\in (1,\frac{p^*}{m})$ such
 that $$\tau M(\tau)\leq \gamma \widehat{M}(\tau),\;\forall \tau\geq \tau_0, \mbox{ where } \widehat{M}(\tau)=\int_{0}^{\tau} M(z)dz.$$
 As in \cite{CP}, we covered the degenerate case. Precisely, we assume:\\\\
$(M_2):$  for each  $\eta>0 \mbox{ there is } m_\eta>0 \mbox{ such that }
 M(\tau)\geq  m_\eta, \forall \tau\geq\eta.$\\

 As $K \in L^\infty(\O)$ then we may assume that $p>1$.
Set $F(s)=\int^s_0 f(t)$ and $p'=\frac{p}{p-1}$ the conjugate exponent of $p$. We introduce the following assumptions on the nonlinear term:
 \begin{enumerate}
\item[$(H_1):$]  there exists $C>0$ such that $C|f(s)|^{q^*}\leq s f(s)- m\gamma F(s)$,
 $\forall s\in\R;$
 \item[$(H_2):$] $\left\{\begin{array}{ll}
\lim\limits_{s \rightarrow \infty}$ $\frac{f(s)}{| s |^{p^*-1}} = 0 \;&\mbox{ if } \O \mbox{ is a bounded domain} \mbox{ and there exists } C > 0 \mbox{ such that } \forall s \in [-1,1]\\
 |f(s) |\leq \; C|s|^{m-1} \;&\mbox{ if } \O \mbox{ is an unbounded domain and } a>0; \\
 |f(s) |\leq \; C|s|^{p^*-1}\;&\mbox{ if } \O=\R^N \mbox{  and  } a=0.\end{array}\right.$

 \item [$(H_3)$:] $\left\{\begin{array}{ll}
                                 \lim\limits_{s\to
\infty} \dfrac{F( s)}{|s|^{m\gamma } }= \infty\; \mbox{ if } \O \mbox{ is a bounded; or  an unbounded domain and } a>0; \\
                                 \mbox{ there exists } \alpha' \in (1, \min(p',\frac{p^*}{m\gamma})) \mbox{ such that }\liminf\limits_{s\to
\infty} \frac{F(s)}{|s|^\frac{p^*}{\alpha'}}>0,  \mbox{ if } \O=\R^N \mbox{  and  } a=0.\end{array}\right.$
\end{enumerate}
Note that the conditions near $0$ in $(H_2)$ are required to only allow the variational setting of \eqref{10}.

Observe first that if we only assume that $M \in C([0,+\infty))$, $K \in L^\infty(\O)$ and $f \in C(\R)$ satisfying $(H_2)$, then the Euler-Lagrange functional associated to problem \eqref{10} is defined by
\begin{eqnarray}\label{kir func}
I(u)=\frac1m \widehat{M}(|D_r u|^{m} +a|u|^m)-\int_\O K(x)F(u),\;\forall u\in E
\end{eqnarray}
and $I\in C^1(E)$ with
\begin{eqnarray}\label{kir deriv}
\langle I'(u), {v}\rangle =M\left(\int_{\O}|D_r u|^{m} +a|u|^m\right)\int_\O(|\textit{D}_r u|^{m-2}\textit{D}_r u\textit{D}_r v +a|u|^{m-1}v)  -\int_\Omega K(x)f(u)v,\; \forall u,v\in E.
\end{eqnarray}
So, $u \in E$ is a weak solution of \eqref{10}
 if and only if $u $ is a critical point of $I$. Problem \eqref{10} is called nonlocal due to the
presence of the term
$ M\left(\int_{\O}|D_r u|^{m} +a|u|^m\right),$
which implies that the equation in \eqref{10} is no longer a pointwise identity.\\

We begin by verifying the Palais-Samale condition; or the Cerami condition under $(H_1)-(H_2)$ (see Remark \ref{r3} for further comments).

\begin{prop} \label{th13}
Assume that $(H_1)-(H_2)$ and $(M_1)-(M_2)$ hold. If, in addition $K$ satisfies $(K_p)$, then
\begin{itemize}
  \item[$1)$] $I$ satisfies the Palais-Smale condition
if $m\geq 2$;
  \item[$2)$] $I$ satisfies the Cerami condition if $1<m<
2$.
\end{itemize}
\end{prop}
We note that for $m \neq 2$, the variational setting of \eqref{10} lacks an ordered Hilbert
space structure and so provokes some mathematical difficulties to obtain infinitely many solutions when $f(.)$ is an odd function. In a famous paper \cite{CP} Colasuonno-Pucci established  multiplicity results by using minimax approach under a restrictive growth condition at zero (see also \cite{DJM}). Differently to \cite{CP,DJM} we will not here impose any control on $f$ at zero if $\O$ is a bounded domain, however if $\Omega$ is an unbounded domain a standard condition near $0$ is required  to ensure that the functional $I$ is well defined. Also, we weaken the Ambrosetti-Rabinowitz condition assumed in \cite{BLW} by assumption $(H_1)$ (see remark \ref{r3}) and we will not employ any cut-off technique as in \cite{BLW,LLS0,LLS1}. Our proof is a direct application of the following symmetric mountain pass theorem.
\begin{theorem} \label{th45}(\cite{R}).
 Let $E$  be a real infinite dimensional Banach space and $I\in C^1(E)$ satisfying the Palais-Smale condition; or the Cerami condition. Suppose $E=E^-\oplus E^+$, where
 $E^-$ is finite dimensional, and assume the following conditions:
\begin{enumerate}
  \item $I$ is even and $I(0)=0$;
 \item there exist $\alpha >0$ and $\rho>0$ such that $I(u)\geq \alpha$ for any $u\in E^+$ with $\|u\|=\rho$;
\item for any finite dimensional subspace $W\subset E$ there is $R= R(W)$ such that $I(u)\leq 0$ for $u\in W$, $\|u\|\geq R$;
\end{enumerate}
then, $I$ possesses an unbounded sequence of critical values.
\end{theorem}

The novelty here consists in using a Schauder basis of $E$ (see Corollary $3$ in \cite{FJN}) to verify point $2$ under only condition $(H_2)$. Precisely, we have

\begin{lem} \label{lem113}
Assume that $(H_2)$ and $(M_2)$ hold. If, in addition $K$ satisfies $(K_p)$, then, {\bf for every} $\rho>0$ there exist a finite dimensional subspace $E^-$ and $\alpha >0$  such that $I(u)\geq \alpha$ for any $u\in E^+$ with $\|u\|=\rho$, where $E^+$ is a topological complement of $E^-$, i.e, $E=E^-\oplus E^+$.
\end{lem}
Our multiplicity existence result reads as follows
\begin{thm}\label{th49}
Suppose that:
\begin{itemize}
 \item $M$ verifies $(M_1)$ and $(M_2)$,
  \item $f$ is an odd function verifying $(H_1)$-$(H_3)$,
  \item $K\in L^\infty(\O)$ is a positive weight function and $K$ satisfies $(K_p)$ if $\O$ is an unbounded,
\end{itemize}
then, $I$ admits infinitely many distinct
pairs $(u_j, -u_j),\;j \in \N^*$, of critical points. Moreover,
$I(u_j)$ is unbounded.
\end{thm}
\subsection{\bf  Mountain pass solution:}
In the following we assume that $\O$ is a bounded domain and we study problem \eqref{10} in the general case $K(x)f(u)=f(x,u)$, that is
\begin{eqnarray} \label{111}
\begin{cases}
M\left(\int_{\O}|D_r u|^{m} \right)\Delta^r_m u = f(x,u) &\mbox{in}\quad \Omega, \\
u=\left(\frac{\partial}{\partial \nu}\right)^k u=0, \quad &\mbox{on}\quad \partial\Omega, \quad k=1, 2,.....  , r-1,
\end{cases}
\end{eqnarray}
For the sake of simplicity, we denote again the following similar conditions by $(H_1)$, $(H_2)$ and $(H_3)$:
\begin{enumerate}
\item[$(H_1):$]  there exist $s_0>0$ and $C>0$ such that $C|f(x,s)|^{q^*}\leq s f(x,s)- m\gamma F(x,s)$,
 $\forall |s|>s_{0}$ and $x\in \Omega,$\\
 where $F(x,s)=\int^s_0 f(x,t)$;
\item[$(H_2):$] $\lim\limits_{s \rightarrow \infty}$ $\frac{f(x ,s)}{| s |^{p^*-1}} = 0$, uniformly with respect to $x \in \Omega$;
\item[$(H_{3}):$] $\lim\limits_{s\to
\infty} \dfrac{F(x, s)}{|s|^{m\gamma } }= \infty$, uniformly in $ \overline\Omega.$
\end{enumerate}
\begin{rem}\label{r22}
We may see easily that Proposition \ref{th13}, Lemma \ref{lem113} and Theorem \ref{th49} hold under the above assumptions.
\end{rem}

 To provide the mountain pass structure in the more familiar setting in the literature in which $M=1$, we require that $F(x,s)$ grows less rapidly than $ \frac{\lambda_1}{m}|s|^m$ near $0$ and more rapidly than $ \frac{\lambda_1}{m}|s|^m$ at infinity, where \begin{eqnarray}\label{lambda1}
\lambda_1:=\inf_{\substack{u\in E\\u\neq 0}}\dfrac{\int_\O|\textit{D}_ru|^m dx }{\int_\O|u|^m dx}>0,\end{eqnarray}
 is the first "eigenvalue" of $\Delta^r_m$. By analogy with $ \lambda_1$, we set
\begin{eqnarray}\label{lambdaM}\lambda_M :=\inf_{\substack{u\in E\\u\neq 0}}\dfrac{  \widehat{M}(\|u\|^m) }{\int_\O|u|^{m\gamma} dx},\end{eqnarray}
with $1<m\gamma <p^*$. To ensure that $\lambda_M $ is positive, we need the following  {\bf sufficient and necessary } coercivity condition: \\\\
$(M_3):$  there is a positive constant $C$ such that $C\tau^\gamma \leq \widehat{M}(\tau), \; \forall \tau \geq 0.$\\
\begin{lem}\label{eig.posi}
\begin{enumerate}
  \item[$(i)$] $\lambda_M  \mbox{ is positive {\bf if and only if} } M \mbox{ satisfies }(M_3);$\label{lam pos}
  \item[$(ii)$] if $M= C\tau^{\gamma-1}$, then $\lambda_M $ is attained.\label{resul2}
\end{enumerate}
\end{lem}
An instructive example will be given in the end of the proof of Lemma \ref{eig.posi}, where $\lambda_M$ is not attained.

To derive a mountain pass solution we assume that $(M_1)$ is global, i.e. $\tau_0=0$, we also replace $(H_3)$ by the following conditions at infinity and at zero:
\begin{enumerate}
\item[$(H'_3):$] $\limsup \limits_{s\to 0} \dfrac{F(x, s)}{|s|^{m\gamma }}< \dfrac{\lambda_{M}}{m } < \liminf\limits_{s\to
\infty} \dfrac{F(x, s)}{|s|^{m\gamma } }, \;{  uniformly \; in }\; \;\overline\Omega.$
 \end{enumerate}
Then, we have
\begin{thm} \label{th24}
Assume that $(H_1)$-$(H_2)$-$(H'_3)$, $(M_1)$ (with $\tau_0=0$) and $(M_3)$ hold.  Then, problem \eqref{111} has a nontrivial mountain pass solution.
\end{thm}
\begin{rem}
 Theorem \ref{th24} holds if we substitute $(M_3)$ by $(M_2)$ and $(H'_3)$ by the following strong condition:
$$\limsup \limits_{s\to 0} \dfrac{F(x, s)}{|s|^{m\gamma }}=0 \mbox{ and }\liminf\limits_{s\to
\infty} \dfrac{F(x, s)}{|s|^{m\gamma } }=\infty, \;{  uniformly \; in }\; \;\overline\Omega.$$
\end{rem}

We close this section by giving some remarks and instructive examples in order to understand the improvement brought by our assumptions.
\begin{rem}\label{r3}
\begin{enumerate}
\item  When $\O=\R^N$ and $a=0$, we may substitute $(H_3)$  by the following condition $\lim\limits_{s\to
\infty} \dfrac{F( s)}{|s|^{m\gamma } }= \infty$, with $m\gamma> \frac{p^*}{p'}$.
 \item  The two main assumptions that appeared in a rich literature ensuring the $(PS)$ condition  are the analogue of the Ambrosetti-Rabinowitz condition related to the Kirchhoff function $M$ (see \cite{CP, DJM, MMTZ, R}):
\begin{enumerate}
\item[${\bf(AR)_\gamma}:$] there are constants $\theta > m\gamma$ and $s_{0} > 0$ such that
 $ s f(x , s) \geq\theta F(x,s)>0,\; \forall  |s| > s_{0} \mbox { and } \forall x \in \Omega;$
\end{enumerate}
\begin{enumerate}
\item[${\bf(SCP)}:$]  there exist $C > 0$ and $p$ satisfying \;$\theta
 \leqslant p < p^*-1$ such that\; $| f(x , s)|  \leqslant C( | s |^{p } + 1), \; \forall (x ,s) \in \Omega  \times \mathbb{R}.$
\end{enumerate}
When $M\equiv1$ (and so $\gamma=1$), some attempts were made to relax conditions $(AR)_1$ and $(SCP)$ (see  \cite{BLW, CM, CSY, H, LL, LRRZ, LZ1, MS, SZ2} and the references therein). Point out that $(H_1)-(H_2)$ are weaker than $(AR)_\gamma$-$(SCP)$. In fact, $(SCP)$ implies $(H_2)$ and from  $(AR)_\gamma$ we have
 $0<(1-\frac {m\gamma}{\theta})sf(x,s)< sf(x,s) -m\gamma F(x,s),\; \forall |s|>s_0$ and $x\in \O.$ As $p^*-1=\frac{1}{q^*-1}$ then from $(H_2)$ we have   $|f(x,s)|^{q^*}\leq sf(x,s), \mbox{  for  } |s|\geq s'_0.$ Clearly, $(H_1)$ follows from the above inequalities.
\\
Note that, $(AR)_\gamma$ requires the following severe restriction called the strong $m\gamma$-superlinear condition:
\begin{enumerate}
  \item [${\bf(SSL)}$:] there exists $C>0$ such that\;\;
$F(x,s)\geq C |s|^\theta,\;\forall x\in \Omega$ and $\forall |s|\geq s_0.$
\end{enumerate}
The most part of the literature used  conditions $(SCP)$ and $(SSL)$ to verify respectively points $2$ and $3$ of Theorem \ref{th45}  which are here weaken by $(H_2)$ and $(H_3)$.
\item Since we assume $\gamma<\frac{p^*}{m}$, then $(m\gamma-1)q^*\leq m\gamma$, so for $(m\gamma-1)q^*-1\leq\alpha<1$ and  $a > \gamma \lambda_M$, then
       a simple computation shows that $f_1(s)=a|s|^{m\gamma-2}s-|s|^{\alpha -1}s$ satisfies $(H_1)$ (and $(H_2)$-$(H'_3)$) but never $(AR)_\gamma$ nor $(SSL)$.
 \item In \cite{SZ2},  $(AR)_1$ was relaxed  into one of the following conditions (for $m=2$ and $M=1$): there are
constants $\theta > 2$ and $C > 0$ such that
 \begin{align}\label{bcz}
  | \theta F(x,s)-s f(x , s) |\leq C(1+s^2),\;\forall (x,s)\in \O\times \R,
  \end{align}
 or the {\bf global} convexity condition
 \begin{align}\label{bdz}
 H(x,s):= sf(x,s)-2F(x,s) \mbox{ is convex in } s,\;\forall x\in\Omega.
  \end{align} However, the following nonlinearity $f_2(s)=|s|^{m\gamma-2}s\ln^q(|s|)$ with $q\geq 1$ verifies
   $(H_1)-(H_3)$ and $(H'_3)$ but does not satisfies assumptions \eqref{bcz} and \eqref{bdz} if $q>1$.\qed
\item  Let $ \gamma_1\in (1,\frac{p^*}{m}),\; \gamma_2 \geq \frac{p^*}{m}$. Consider the degenerate Kirchhoff function $M(\tau)=\tau^{\gamma_1-1} \mbox{ if } \tau\geq 1$ and
    $M(\tau)=\tau^{\gamma_2-1} \mbox{ if } \tau\leq 1$. We can see that $M$ satisfies $(M_1)-(M_2)$ but not the global assumption $(M)$ required in \cite{CP}.
\item
Consider the following Kirchhoff function introduced in \cite{CP}:
  $$ M(\tau)= a\tau^{\gamma_1-1} +b \tau^{\gamma_2-1} \mbox{ with } a\geq 0, b > 0 \mbox{ and } 1\leq \gamma_1\leq \gamma_2 <\frac{p^*}{m}.$$
Then $M$ satisfies $(M_1)-(M_2)$ and also assumption $(M_3)$. Moreover $M$ is degenerate if $a=0$ or $a\neq0$ and $\gamma_1\neq 1$.

 \end{enumerate}
\end{rem}

The outline of this paper is the following: In Section \ref{section 2}, we give proofs of Proposition \ref{th13}, Lemma \ref{lem113} and Theorem \ref{th49}. Section \ref{section 3} is devoted to the proofs of Lemma \ref{eig.posi} and Theorem \ref{th24}.

In the following, $|\cdot|$ denotes the Lebesgue measure in $\R^N$ and  $C$ (respectively $C_{\epsilon}$) denotes always a generic positive constant independent of $n$ and $\epsilon$ (respectively independent of $n$), even their value could be changed from one line to another one.

\section{Proofs of Proposition \ref{th13}, Lemma \ref{lem113} and Theorem \ref{th49} \label{section 2}}
\subsection{ \bf Preliminary results}
Let-us first establish some inequalities from assumptions $(M_1)$ and $(H_1)$-$(H_2)$ which will be useful to prove Proposition \ref{th13} and Lemma \ref{lem113}.

Recall that as $K \in L^\infty(\O)$ then we may assume in assumption $K_p$ that $p>1$. Denote $p'=\frac {p}{p-1}>1$  the conjugate exponent of $p$.
 According to $(M_1)$ and $(H_1)$ , there exists $C_0>0$ such that
\begin{eqnarray}\label{er}
C|f(s)|^{q^*} \leq s f(s)- m\gamma F(s),\; \forall s\in  \R, \mbox{ and }
-C_0 \leq\gamma \widehat{M}(\tau)-\tau M(\tau), \; \forall \tau\geq 0.
\end{eqnarray}
  From $(H_2)$ it follows that for any $\epsilon>0$ there exists $C_{\epsilon}>0$ such that
\begin{eqnarray}\label{k10}
|f(s)|\leq\begin{cases}
  \epsilon|s|^{p^*-1}+C_{\epsilon}, & \mbox{if } \O \mbox{ is a bounded domain};\\
  \epsilon|s|^{p^*-1}+C_{\epsilon}|s|^{\frac{m-1}{p'}}, & \mbox{if } \O \mbox{ is an unbounded domain};\\
   \epsilon|s|^{p^*-1}+C_{\epsilon}|s|^{\frac{p^*-1}{p'}}, & \mbox{if }\O=\R^N \mbox{ and } a=0;
\end{cases}
\end{eqnarray}
\begin{eqnarray}\label{pp}
|F(s)|\leq\begin{cases}
  \epsilon|s|^{p^*}+C_{\epsilon}, & \mbox{if } \O \mbox{ is a bounded domain};\\
  \epsilon|s|^{p^*}+C_{\epsilon}|s|^{\frac{m}{p'}}, & \mbox{if } \O \mbox{ is an unbounded domain};\\
   \epsilon|s|^{p^*}+C_{\epsilon}|s|^{\frac{p^*}{p'}}, & \mbox{if }\O=\R^N \mbox{ and } a=0.
\end{cases}
\end{eqnarray}
 Multiplying \eqref{k10} by $|s-s'|, s'\in \R$, and applying Young's inequality, then we derive
 \begin{eqnarray}\label{llvw}
|f(s) (s-s')|\leq\begin{cases}
  \epsilon(|s|^{p^*}+ |s-s'|^{p^*}) +C_{\epsilon},\;\forall s\in \R & \mbox{if } \O \mbox{ is a bounded domain};\\
  \epsilon(|s|^{p^*}+ |s-s'|^{p^*}+|s-s'|^{m}) +C_{\epsilon}|s|^{\frac{m}{p'}},\;\forall s\in \R & \mbox{if } \O \mbox{ is an unbounded domain};\\
   \epsilon(|s|^{p^*}+ |s-s'|^{p^*}) +C_{\epsilon}|s|^{\frac{p^*}{p'}},\;\forall s\in \R & \mbox{if }\O=\R^N \mbox{ and } a=0;
\end{cases}
\end{eqnarray}
We also recall some known results which will be essential to prove Proposition \ref{th13}. Consider the functional $\psi(u)= \frac{1}{m}\|u\|^{m},\; u\in E$, we have $ \psi \in C^1(E)$ with Fr\'echet's derivative
$$\langle \psi'(u), v\rangle=\int_{\R^N}|\textit{D}_r u|^{m-2}\textit{D}_r u\textit{D}_r v +a|u|^{m-2}uv, \forall u,v \in E$$
where $\O$ is a bounded domain and  the zero mass case $\O=\R^N$.\\
Set $\varphi(t)=t^{m-1}, t\geq 0$, clearly we have $\langle \psi'(u), u\rangle=\varphi(\|u\|)\|u\|$ and  it follows from H\"older's inequality that $\|\psi'(u)\|_{E^*}=\varphi(\|u\|)$.
Obviously, $\varphi$ is a normalization function and since $E$ is locally uniformly convex and so uniformly convex and reflexive Banach space, then the corresponding duality mapping $J_\varphi$ is single valued (i.e., $J_\varphi=\psi'$) and satisfies the $S_+$ condition (see Proposition $2$ in \cite{DJM}, respectively Lemma $3.2$ in \cite{LZ0}):
\begin{eqnarray}\label{khc}
\mbox{ if } u_n\rightharpoonup u \mbox{ and }\limsup_{n\to +\infty}\psi'(u_n)(u_n-u)\leq 0,\;\mbox{ then } u_n \rightarrow u.
\end{eqnarray}

\subsection{ \bf Proof of Proposition \ref{th13}.}
 Since we assume that $N> rm$, we may easily see that
 \begin{eqnarray}\label{err}
 m > \frac{1}{q^*}+1 \mbox{ if  } m\geq 2.
\end{eqnarray}
Let $u_n$ be a $(PS)$ sequence of $I$ if $m\geq 2$ ({ respectively} $(C)$ sequence if $1< m< 2$). The case $u_n$ admits a subsequence which converges strongly to $0$ in $E$ is trivial. Hence, we may suppose that there exist $\eta_0>0$ and $n_0 \in \N$ such that $\|u_n\|^m\geq \eta_0$ for all $n\geq n_0$. So, in view of $(M_2)$ we can find $m_{\eta_0}>0$ such that
 \begin{eqnarray}\label{ing1}
 m_{\eta_0}\|u_n\|^m \leq M(\|u_n\|^m)\|u_n\|^m, \forall n\geq n_0.
\end{eqnarray}
\\
\textbf{Step 1.} We shall prove that $u_n$ is bounded in $E$. First, from \eqref{kir deriv} and  \eqref{ing1}, we have
\begin{align*}
m_{\eta_0}\|u_n\|^m \leq  &\langle I'(u_n),u_n\rangle+\int_{\O} Kf(u_n)u_n.
\end{align*}
Apply H\"older's inequality to the second term in the right-hand side and using \eqref{Gagliardo}, we obtain
\begin{align}\label{kirch hk}
 m_{\eta_0}\|u_n\|^m \leq \langle I'(u_n),u_n\rangle&+ C\left(\int_{\O} K|f(u_n)|^{q^*}\right)^\frac{1}{q^*}\|u_n\|.\end{align}
From \eqref{kir func} and \eqref{kir deriv}, one has
\begin{eqnarray}\label{ing2}
\int_\O K \left [ f(u_n)u_n-m\gamma F(u_n)\right]= \left[\gamma
\widehat{M}(\|u_n\|^m)-M(\|u_n\|^m)\|u_n\|^m\right]+ \langle I'(u_n), {u_n}\rangle-m\gamma I(u_n).
 \end{eqnarray}
Taking into account that $u_n$ is a $(PS)$ sequence if $m\geq 2$ ( respectively $(C)$ sequence if $1< m< 2$), then from \eqref{er} we deduce
\begin{equation*}\int_{\Omega} K| f( u_{n}) |^{q^*}\leq C(1+\|u_{n}\| ), (\mbox{  respectively }\int_{\Omega}K | f( u_{n})
|^{q^*}\leq  C).
  \end{equation*}
 Combining now the above inequality with \eqref{kirch hk}, it follows
\begin{equation*}
\|u_{n}\|^{m}\leq C( 1 + || u_{n} ||^{{\frac{1}{q^*}+1}} )\mbox{ if } m\geq 2, (\mbox{  respectively } \|u_{n}\|^{m}
\leq C \| u_{n} \| \; \mbox{ if } 1<m<2).
\end{equation*}
 Therefore, thanks to \eqref{err} the $(PS)$ sequence is bounded in $E$ if $m\geq 2$, and clearly, the $(C)$ sequence $u_n$ is also bounded if $1<m<2$.
 \\\\
{\bf Step 2.} We shall prove that the bounded sequence $u_n$ has
a strong convergent subsequence in $E$. In fact, we can find a subsequence (denoted again by $u_n$) and $u \in E$ such that $u_n$ converges to
$u$ weakly in $E$, \text{ and a.e. in } $\O$. Also $u_n$ and $u_n - u$ are bounded in $L^{p^*}(\O)$ (and in $L^{m}(\O)$ if $\O$ is an unbounded domain and $a>0$). Using again \eqref{kir deriv}, we get
\begin{eqnarray}\label{e:1221zz}
M(\|u_n\|^m)\int_\O\left(|\textit{D}_r u_n|^{m-2}\textit{D}_r u_n\textit{D}_r (u_n-u) +a|u_n|^{m-2}u(u_n-u)\right) =\langle I'(u_n), u_n-u\rangle+\int_\Omega K f(u_n)(u_n-u).
\end{eqnarray}
First we claim that
 \begin{eqnarray}\label{ggf}
                     \int_\Omega Kf(u_n)(u_n-u)\mbox{ converges to } 0.
                     \end{eqnarray}
 To prove \eqref{ggf} we shall discuss three cases:\\
{\bf Case 1. $\O$ is a bounded domain.} We may assume that $u_n$ converges to $u$ in $L^1(\O)$.
Apply inequality \eqref{llvw} (with $s= u_n$ and $s'=u$, and integrate over $\O$, we obtain
 $$\left|\int_\Omega K f(u_n)(u_n-u)  \right|\leq \epsilon + C_\epsilon \int_{\O}|u_n-u|.$$
As $u_n$ converges strongly in $L^1(\O)$, then \eqref{ggf} follows.\\
{\bf Case 2. $\O$ is an unbounded domain and $a>0$.} From inequality \eqref{llvw}, we have
\begin{eqnarray}\label{fxfR}
K(x)|f(u_n)(u_n-u)|\leq \epsilon K_\infty (|u_n|^{p^*}+|u_n-u|^{p^*}+|u_n-u|^{m})+ C_{\epsilon} K|u_n|^{\frac{m}{p'}} & \mbox{if } \O \mbox{ is an unbounded domain}.
\end{eqnarray}
Let $A\subset \R^n$ be measurable set such that $|A|^{\frac 1 p}<\frac{\epsilon}{ K_{\infty}C_{\epsilon}}$ where $K_{\infty}=\|K\|_{L^{\infty}(\O)}$.  Integrate \eqref{fxfR} over $A$ and using H\"older's inequality, we deduce
\begin{eqnarray}\label{ffA}
\int_{A}K(x) |f(u_n)(u_n-u)|
 \leq C \epsilon  +  C_{\epsilon} K_{\infty} |A|^{\frac 1 p}\left(\int_{A}|u_n|^{m}\right)^{ \frac{1}{p'}}\leq C'\epsilon.
\end{eqnarray}
Invoking now the assumption $(K_p)$, then we may choose $R_\epsilon>0$ large enough such that $$\left(\int_{\{|x|> R_\epsilon\}}K^{p}(x)\right) ^{\frac 1 p} \leq \frac{\epsilon}{C_{\epsilon}}.$$ Integrating \eqref{fxfR} over $\{|x|> R_\epsilon\}$, we derive
\begin{eqnarray}\label{ffR}
\int_{\{|x|> R_\epsilon\}}K(x)|f(u_n)(u_n-u)| \leq C\epsilon  + C_{\epsilon} \int_{\{|x|> R_\epsilon\}}K(x)|u_n|^{ \frac{m}{p'}}
 \leq C\epsilon  + C_{\epsilon} \left(\int_{\{|x|> R_\epsilon\}}K^{p}(x)\right) ^{\frac 1 p}\left(\int_{\{|x|> R_\epsilon\}}|u_n|^{m}\right)^{ \frac{1}{p'}}\leq C'\epsilon.\nonumber\\
\end{eqnarray}
 Also, $K(x)f(u_n)(u_n-u) \to 0 \text{ a.e. in } \R^N$, then in view of  \eqref{ffA} and \eqref{ffR}  and Vitali's theorem we deduce \eqref{ggf}.\\
{\bf Case 3. The $m\gamma-$zero mass $\O=\R^N$ and $a=0$.} Recall first \eqref{llvw} in this case by
\begin{eqnarray*}
|f(s) (s-s')|\leq\epsilon(|s|^{p^*}+ |s-s'|^{p^*}) +C_{\epsilon}|s|^{\frac{p^*}{p'}},\;\forall s\in \R,
\end{eqnarray*}
as above we deduce  $$K(x)|f(u_n)(u_n-u)|\leq \epsilon K_\infty (|u_n|^{p^*}+|u_n-u|^{p^*})+ C_{\epsilon} K|u_n|^{\frac{p^*}{p'}}.$$
Similarly to the case $2$ (where we have  only to substitute $m$ by $p^*$ and $a$ by $0$) we deduce \eqref{ggf}.

Next, since $I'(u_n) \to 0$ in $E^*$ and $(u_n-u)$ is bounded in $E$, then from \eqref{e:1221zz} and \eqref{ggf} we can deduce that $$M(\|u_n\|^m)\int_\O\left(|\textit{D}_r u_n|^{m-2}\textit{D}_r u_n\textit{D}_r (u_n-u) +a|u_n|^{m-2}u(u_n-u)\right) \to 0.$$
Hence,  \eqref{ing1} yields  $\int_\O\left(|\textit{D}_r u_n|^{m-2}\textit{D}_r u_n\textit{D}_r (u_n-u) +a|u_n|^{m-2}u(u_n-u)\right)  \mbox{ converges to } 0.$ Invoking now the $S_+$ property (see \eqref{khc}), we conclude that $u_n$ converges strongly to $u$ in $E$. \qed
\subsection{ \bf  Proof of Lemma \ref{lem113}.}
 Let $(e_i)_{i \in \N^*}$ be a Schauder basis of $E$ (see Corollary $3$ in \cite{FJN} and also \cite{FHHSPZ,T}), which means that each $x\in E$ has a {\bf
unique representation} $x=\sum_{i=1}^{\infty}a_ie_i$, where $a_i$ are real numbers. Set $ E_j=span(e_1, e_2,..,e_j)$, then, the linear projection onto $ E_j$ i.e.,
$P_j: E\to E_j,\; P_j(x)= \sum_{i=1}^{j}a_ie_i$ is a {\bf continuous} linear operator for all $j\in \N^*$ (see \cite{S2}) \footnote{More precisely, $P_j$ are
 uniformly bounded, that is there exists $C>0$ such that $\|P_j(x)\|\leq C\|x\| $ for each $j\in \N^*$ and all $x\in E$ (see \cite{FHHSPZ,S2})}.
  Therefore, $F_j= N( P_j)$ (the kernel of $P_j$) is a topological complement of $E_j$, that is $E=E_j\oplus F_j$.\\
    Fix $\rho\geq 0$ and set $S_j^\perp(\rho)=\left\{u\in F_j \mbox{ such that }\|u\|=\rho\right\}$ and $\beta_j:=\sup_{S_j^\perp(\rho)} \int_{\Omega} K|F(u)|$.\\
  So, we claim that
 \begin{eqnarray}\label{uxy}
 \beta_j \rightarrow 0, \mbox{ as } j\rightarrow +\infty,
 \end{eqnarray}
  Before proving the claim \eqref{uxy}, let-us first end the proof of Lemma \ref{lem113}. Fix $\rho >0$, then for all $u\in S_j^\perp(\rho)=\left\{u\in F_j \mbox{ such that }\|u\|=\rho\right\}$, we have
$$I(u)\geq \frac 1m \widehat{M}( ||u||^m)-\int_\O K F(u)\geq \frac 1m \widehat{M}( \rho^m) -\beta_j.$$
Set $\alpha= \frac {1}{2m} \widehat{M}(  \rho^m)$, since we assume that $M$ satisfies $(M_2)$, then $\alpha >0$.
As $\beta_j$ converges to $0$, we can choose $j=j_0$ large enough such that $\beta_{j_0} \leq  \alpha $. Hence,
 $I(u)\geq  \alpha$, and so Lemma \ref{lem113} holds with $E^-=E_{j_0}$, $E^+=F_{j_0}$ and $\alpha= \frac {1}{2m} \widehat{M}(  \rho^m)$.
\\\\
{\bf Proof of \eqref{uxy}.}
 We argue by contradiction. Suppose that there exist $m_0>0$ and a subsequence (denoted by $\beta_j$) such that\\
$$m_0<\beta_j,\forall j\in \N^*.$$
From the definition of $\beta_j$, there exists $u_j\in S_j^\perp(\rho)$ such that
\begin{eqnarray}\label{uuk}
m_0<\int_{\O} K|F(u_j)|\leq \beta_j.
\end{eqnarray}
   As $\|u_j\|=\rho$, then there exist a subsequence (denoted by $u_{j}$) and $u\in E $ such that $u_{j}$  converges weakly to $u$ and $\mbox{ a.e in }\; \Omega$, also $u_j$ is bounded in $L^{p^*}(\O)$. Fix $k\in \N^*$, as $P_k\circ P_{k+1}=P_k$, then $ F_{j}\subset F_k$ for all $j\geq k$, and so
 \begin{eqnarray}\label{uuxk}
u_{j} \in F_{k},\mbox{ that is } P_k(u_{j})\; \forall j\geq k.
\end{eqnarray}

 On the other hand, Since $P_k$ is a continuous linear operator and $u_{j}$  converges weakly to $u$, then $P_k(u_{j})$
 converges weakly to $P_k(u)$ which with \eqref{uuxk} implies that $P_k(u)=0$ and therefore $u=0$ as $u=\lim_{k\to
\infty}\sum_{i=1}^{k}a_ie_i=\lim_{k\to
\infty}P_k(u)$. Consequently, $u_j$ converges to $0$ $\mbox{ a.e in } \Omega$, and then $KF(u_{j})  \rightarrow 0 \mbox{ a.e in }\; \Omega$ as $KF(0)=0$. Substitute $s$ by $u_j$ in inequality \eqref{pp} and integrate over a measurable set $A\subset\O$, then if $\O$ is a bounded domain we derive
\begin{eqnarray*}
\int _{A} K|F(u_{j}) | \leq   C\epsilon  \int _{A} | u_{j}|^{p^*}
+ C_{\epsilon}|A| \leq C\epsilon.
\end{eqnarray*}
Hence if  $|A|<\frac{\epsilon}{C_{\epsilon}}$ (where $|A|$ denotes the Lebesgue measure of $A$), we deduce
\begin{eqnarray*}
\int _{A} K|F(u_{j}) |  \leq C'\epsilon.
\end{eqnarray*}
Taking into account that $\O$ is a bounded domain, then Vitali's theorem implies that $K|F(u_{j})|\rightarrow 0$ in $L^{1}(\O)$ and in view of \eqref{uuk}, we obtain
$$
0< m_0\leq 0.
$$
 Thus, we reach a contradiction.  Now, we consider the case of unbounded domains and $a>0$, and in which assumption $(H_2)$ implies:\\
$$|F(s)|\leq\epsilon|s|^{p^*}+C_{\epsilon}|s|^{\frac{m}{p'}},$$
which shows that
\begin{align*}
\int _{A} |KF(u_{j})| &\leq   \epsilon K_\infty \int _{A} | u_{j}|^{p^*}+
+ K_\infty C_{\epsilon}\int _{A} |u_{j}|^\frac{m}{p'}\\
&\leq C\epsilon+ K_\infty C_{\epsilon} |A|^{\frac 1 p} \left(\int _A|u_{j}|^{m}\right)^\frac{1}{p'}
\\&\leq C(\epsilon+  C_{\epsilon} |A|^{\frac 1 p}).
\end{align*}
Therefore, if $|A|^{\frac 1 p}<\frac{\epsilon}{K_\infty C_{\epsilon}}$, we deduce
\begin{eqnarray}\label{uxuW}
\int _{A} K|F(u_{j}) |  \leq C'\epsilon.
\end{eqnarray}
By the virtue of $(K_p)$  there is $R_\epsilon>0$ such that $\left(\int_{\{|x|> R_\epsilon\}}K^{p}(x)\right) ^{\frac 1 p} \leq \frac{\epsilon}{C_{\epsilon}}$. As above, we have
\begin{align*}
\int _{\{|x|> R_\epsilon\}} |KF(u_j)| &\leq  \epsilon K_\infty  \int _{\{|x|> R_\epsilon\}} | u_ j|^{p^*}
+ C_{\epsilon} \int _{\{|x|> R_\epsilon\}}K(x)|u_{j}|^\frac{m}{p'} \\
&\leq C\epsilon + C_{\epsilon} \left(\int_{\{|x|> R_\epsilon\}}K^{p}(x)\right) ^{\frac 1 p}\left(\int_{\{|x|> R_\epsilon\}}|u_j|^{m}\right)^{ \frac{1}{p'}} \\&\leq C\epsilon.
\end{align*}
Apply now Vitali's theorem to deduce that $K|F(u_{j})|\rightarrow 0$ in $L^{1}(\O)$ and as above we reach a contradiction from \eqref{uuk}.

For the $m\gamma$-zero mass we omit here the proof which is similar to one of the unbounded domain where we substitute $m$ by $p^*$ (see assumption $(H_2)$). Thus, the proof \eqref{uxy} is well completed.\qed

\subsection{ \bf Proof of Theorems \ref{th49}}
We will show that the functional $I$ satisfies all conditions of the abstract Theorem \ref{th45}. In fact, Since $f$ is odd and $F(0)=0$, then $I$ is an
even functional and $I(0)=0$, and according to Proposition \ref{th13}, $I$ satisfies the $(PS)$ condition if $m\geq 2$ (respectively the $(C)$ condition if $1<m<2$). Thanks to Lemma \ref{lem113} $I$ verifies point $2$ of Theorem \ref{th45}. Therefore, it remains to show that condition $3$ of Theorem \ref{th45} holds. We shall first consider the delicate case:\\\\
{\bf The case $\mathbf{\O=\R^N}$ and $\mathbf{a=0}$.}
As $K \in L^\infty(\R^)$, we may assume in this case that $(K_p)$ is verified for $\mathbf{p\geq q^*}$.  Recall the assumption $(H_3)$ :
\\
$(H_3)$: There exists  $\alpha' \in (1, \min(p',\frac{p^*}{m\gamma})) \mbox{ such that }\liminf\limits_{s\to
\infty} \frac{F(s)}{|s|^\frac{p^*}{\alpha'}}=L>0.$
\\
Let $\alpha$  be the conjugate exponent of $\alpha'$. Since $\alpha'<p'$, then $K\in  L^{\alpha}(\R^N)$ and $\alpha>p>q^*$. So $\frac{p^*}{\alpha'}> 1$
and from H\"older's inequality we have $|u|_{K,\alpha'}=\left(\int_{\R^N}K(x)|u|^{\frac{p^*}{\alpha'}}\right)^\frac{\alpha'}{p^*}$ is well defined. Consequently, for $u,v\in E$ we have $K^\frac{\alpha'}{p^*}u,\; K^\frac{\alpha'}{p^*}v\in L^\frac{p^*}{\alpha'}(\R^N)$, and therefore
\begin{eqnarray*}
 |u+v|_{K,\alpha'}&=&\left(\int_{\R^N}\Big(K(x)^\frac{\alpha'}{p^*}|u+v|\Big)^{\frac{p^*}{\alpha'}}\right)^\frac{\alpha'}{p^*}
 \\&\leq& \left(\int_{\R^N}\Big(K(x)^\frac{\alpha'}{p^*}|u|+K(x)^\frac{\alpha'}{p^*}|v|\Big)^{\frac{p^*}{\alpha'}}\right)^\frac{\alpha'}{p^*}.
 \end{eqnarray*}
 Using the Triangle inequality  in $L^\frac{p^*}{\alpha'}(\R^N)$, we derive
 \begin{eqnarray*}|u+v|_{K,\alpha'}\leq |u|_{K,\alpha'}+|v|_{K,\alpha'}.
\end{eqnarray*}
Consequently, $|u|_{K,\alpha'}$ define a norm in $E$. Now, according to $(H_3)$, there exists $s_L> 0$ such that\\ $F(s)\geq \frac{L}{2}|s|^\frac{p^*}{\alpha'}, \forall |s|>s_L,$ and $F(s)\geq -C|s|^\frac{p^*}{p'}, \forall\; |s| \leq 1.$ Hence, we can find $C_L>0$ such that
\begin{eqnarray*}\label{dww}
F(s)\geq \frac{L}{2}|s|^\frac{p^*}{\alpha'}-C_L|s|^\frac{p^*}{p'},\; \forall s\in \R.\end{eqnarray*}
 In view of $(M_1)$,  we can find $C_1>0$ such that
\begin{eqnarray}\label{cond-Kir}
\widehat{M}(\tau)\leq C_1 \tau^\gamma \mbox{ for all } \tau\geq 1.
\end{eqnarray}
Combine the above inequalities, then for all $\|u\|\geq 1$ we have
\begin{equation*}
I(u)\leq  \frac{C_1}{m}\|u\|^{m\gamma} - \frac{L}{2} |u|_{K,\alpha'}^\frac{p^*}{\alpha'}+C|K|_{L^p(\R^N)}\|u\|^\frac{p^*}{p'}_{L^{p^*}(\R^N)}.
\end{equation*}
Let $W$ be a fixed finite dimensional subspace of $ E$, as $||.||$,
$|u|_{K,\alpha'}$ and $\|u\|_{L^{p^*}(\R^N)}$ are equivalent norms on $W$, we can find $C_W>0$ such that
 $$I(u)\leq  \frac{C_1}{m}\|u\|^{m\gamma} - \frac{L}{2}C_W \|u\|^\frac{p^*}{\alpha'}+\frac{C}{C_W}|K|_{L^p(\R^N)}\|u\|^\frac{p^*}{p'}.$$
Since $\frac{p^*}{\alpha'}>\frac{p^*}{p'}>m\gamma$ then we may find $R= R(W)>1$ large enough such that $I(u)< 0$ for all $\|u\|\geq R, \; u\in W$.
\\
$\\$
{\bf The case $\mathbf{\O}$ is an unbounded domain with $\mathbf{a>0}$.}
\\
Recall that if $a>0$ the functional space $E= W^{r,m}(\O)\hookrightarrow L^s(\O)$ for all $m\leq s\leq p^*$. As $K\in L^{\infty}(\O)$, then $\left(\int_{\O}K(x)|u|^{s}\right)^\frac{1}{s}$ defines a norm. From $(H_3)$, we verify as above that for every $A>0$ there is $C_A>0$ such that
\begin{eqnarray}\label{DW}
F(s)\geq A|s|^{m\gamma}-C_A|s|^{m},\; \forall s\in \R,
\end{eqnarray}
and so for all $\|u\|\geq 1$ we have
\begin{equation*}
I(u)\leq  \frac{C_1}{m}\|u\|^{m\gamma} - A \int_{\O}K(x)|u|^{m\gamma}+C_AK_\infty \int_{\O}|u|^m.
\end{equation*}
Let $W$ be a finite dimensional subspace of $ E$. Taking into account that $\left(\int_{\O}K(x)|u|^{m\gamma}\right)^\frac{1}{m\gamma}$ is a norm and since all norms on $W$ are equivalent, we can find $C_W>0$ such that\begin{equation*}
I(u)\leq  (\frac{C_1}{m} -A C_W) \|u\|^{m\gamma}+ \frac{K_{ \infty }C_A}{ C_W }\|u\|^{m}.
\end{equation*}
Choosing  $A=\dfrac{2C_1}{m C_W}$, as $m\gamma>m$ then we may find $R= R(W)>1$ large enough such that $I(u)< 0$ for all $\|u\|\geq R, \; u\in W$.
\\
{\bf The case $\mathbf{\O}$ is a bounded domain.}
\\\\
In view of $(H_3)$ we have for all $A>0$, there is $C_A>0$ such that $$F(x,s)\geq A|s|^{m\gamma}- C_A,\; \forall (x,s)\in
\O\times \mathbb{R}.$$ Therefore we may conclude $$ I(u) \leq
\frac{C_1}{m}||u||^{m\gamma}-A\|u\|^{m\gamma}_{L^{m\gamma}(\O)} +C_A,\;\forall\|u\|>1.$$
The desired result followed as in the previous case which completes the proof of Theorem \ref{th49}.\qed

\section{Proofs of  Lemma \ref{eig.posi} and Theorem \ref{th24}} \label{section 3}
\subsection{ \bf Proof of Lemma \ref{eig.posi}} { \bf Proof of $(i)$:} Assume that $M$ satisfies $(M_3)$. As $1<m\gamma <p^*$, Sobolev's inequality implies
$\int_\O|u|^{m\gamma} \leq C \|u\|^{m\gamma},$ which combined with $(M_3)$ yields  $\lambda_M >0$.

 Conversely, if $M$ does not satisfies $(M_3)$, then there is a sequence $\tau_i > 0$ such that $\tau_i^{-\gamma}\widehat M(\tau_i) \to 0$. Consider $\varphi \in E$ such that $\|\varphi\| = 1$, set $u_i = \tau_i^{1/m}\varphi$, then $\|u_i\|^m = \tau_i$. So
$$\frac{\widehat M(\|u_i\|^m)}{\int_\O|u_i|^{m\gamma}} = \frac{\widehat M(\tau_i)}{\tau_i^\gamma} \times \frac{1}{\int_\O|\varphi|^{m\gamma} } \to 0.$$ Therefore, $\lambda_M = 0$.\qed\\\\
{ \bf Proof of $(ii)$:} $(a)$ We will first prove that if $M= C\tau^{\gamma-1}$, then $\lambda_M$ is attained. without losing
any generality, we may assume that $C=\gamma$.
Let
\begin{eqnarray*}\label{lambdaM}\lambda_M :=\inf_{\substack{u\in E\\u\neq 0}}\dfrac{  \widehat{M}(\|u\|^m) }{\int_\O|u|^{m\gamma} }=\inf_{\substack{u\in E\\u\neq 0}}\dfrac{ \frac{C}{\gamma}\|u\|^{m\gamma}}{\|u|_{L^{m\gamma}(\O)}^{m\gamma} }=\inf\left\{\|u\|^{m\gamma}, \mbox{ such that } \|u|_{L^{m\gamma}(\O)}=1\right\}.\end{eqnarray*}

Let $u_n$ be a minimizing sequence, i.e, \; $\|u_n\|_{L^{m\gamma}(\O)}=1$ and \;$\|u_n\|^{m\gamma}\to \lambda_M$, so $u_n$ is bounded in the $E$ norm. Therefore, as $1<m\gamma<p^*$, there is $u\in E$ and a subsequence
(still denoted by $u_n$) such that $u_n$ converges weakly to $u$ in $E$, $\|u_n\|_{L^{m\gamma}(\O)}\to \|u\|_{L^{m\gamma}(\O)}$, and
\begin{equation*}
\lambda_M=\liminf_{n\to +\infty} \|u\|^{m\gamma}\geq \|u\|^{m\gamma}.
\end{equation*}
 Consequently, $\|u\|_{L^{m\gamma}(\O)}=1$ and so $\|u\|^{m\gamma}\geq\lambda_M$ which implies that $\|u\|^{m\gamma}=\lambda_M$.
\\
Moreover, there exists a Lagrange multiplier $\mu$ such that

$$mM(\|u\|^m)\int_\O|\textit{D}_r u|^{m-2}\textit{D}_r u\textit{D}_r v  =\mu m\gamma\int_\Omega|u|^{m\gamma-1} v,\; \forall v\in E.$$
If $v=u$ we have $m\lambda_M=\mu m\gamma$, that is $\lambda_M=\mu m\gamma$ and so
\begin{eqnarray*}
\begin{cases}
M(\|u\|^m)\Delta^r_m u =\lambda_M|u|^{m\gamma-1}  &\mbox{in}\quad \Omega, \\
u=\left(\frac{\partial}{\partial \nu}\right)^k u=0, \quad &\mbox{on}\quad \partial\Omega, \quad k=1, 2,.....  , r-1.
\end{cases}
\end{eqnarray*}
$(b)$ In general $\lambda_M$ is not attained if $M$ satisfies $(M_3)$. The typical example is $\widehat{M}(t)=t^\gamma H(t)$, where
$H:\R_+\to \R_+$ is strictly decreasing
with $\lim_{t\to \infty}H(t)=l>0$
  and $t^\gamma H(t)$ is strictly increasing on $\R_+$. So, $M$ is positive and satisfies $(M_3)$.
\\
Set
 $$E(u)=\dfrac{  \widehat{M}(\|u\|^m) }{\int_\O|u|^{m\gamma} dx}.$$
The monotonicity of $H$ involving $E(\alpha u)<E(u)$ \;for all\; $\alpha>1,$ and $u\neq0.$ It means
clearly that $\lambda_M$ is not attained. More exactly, one can
find easily examples of $H$ such as $H(t)=1+\beta(t+1)^{-1}$ (with $\beta>0$) and which also satisfy $(M_1)$.\qed
\subsection{ \bf Proof of  Theorem \ref{th24}}
First of all observe that $(M_1)$ (with $\tau_0=0$) implies
that for each $\tau_1>0$, we have
\begin{eqnarray}\label{MMM}
 \frac{\widehat{M}(\tau)} {\tau^\gamma} \leq \frac{ \widehat{M}(\tau_1)}{\tau_1^\gamma} , \; \forall \tau \geq \tau_1.
\end{eqnarray}
  To prove  Theorem \ref{th24}, we shall verify the validity of the conditions of the standard mountain pass theorem \cite{R}. Since $(M_3)$ implies $(M_2)$, Proposition \ref{th13} holds. Consequently,  $I$ satisfies the $(PS)$ condition if $m\geq 2$ (respectively the $(C)$ condition if $1<m<2$). By combining $(H_{2})$ and $(H'_{3})$ (at $0$), we can find $\epsilon_{0} > 0$ small enough and $C_0>0$ such that
$F(x ,s)\leq(\frac{\lambda_{M} }{m}-\epsilon_{0}) | s |^{m\gamma}+ C_0 | s |^{p^*} \mbox{ for all } (x,s)\in \Omega \times\mathbb{R}.$
 Also recall that $(M_3)$ implies $(i)$ of Lemma \ref{eig.posi} which with \eqref{Gagliardo} implies
 \begin{align*}
 I(u)  &\geq  \frac{1}{m}\widehat{M}(\|u\|^m)  -  (\frac{\lambda_{M} }{m}-\epsilon_{0}) \int_\Omega |u|^{m\gamma} - C_0 \int_\Omega |u|^{p^*}\\&\geq \frac{\epsilon_{0}}{\lambda_{M}}\widehat{M}(\|u\|^m)  - C'_0 || u ||^{p^*},\; C'_0>0.
 \end{align*}
Set $ \|u\|=\rho$ with $0<\rho \leq 1$, thus using $(M_3)$, we deduce
\begin{align*}
 I(u)  \geq \frac{C\epsilon_{0}}{\lambda_{M}}\rho^{m\gamma} - C'_0 \rho^{p^*}\geq \rho^{m\gamma}(\frac{C\epsilon_{0}}{\lambda_{M}} - C'_0 \rho^{p^*-m\gamma}).
 \end{align*}
     Choose $\rho =\inf(1,(\frac{C\epsilon_{0}}{2C'_0\lambda_{M}})^{\frac{1}{p^*-m\gamma}})$ and $\alpha=\frac{C\epsilon_{0}}{2\lambda_{M}}\rho^{m\gamma} >0$, then we have $I(u)\geq \alpha$ for all $\|u\|=\rho$.

 On the other hand, using $(H'_{3})$ (at infinity) and part (i) of Lemma \ref{lam pos}, then for $\epsilon_0>0$ small enough, we can find a positive constant $C_0$ and $\varphi \in E\setminus\{0\}$ such that
\begin{eqnarray}\label{e44}
|F(x,s)| \geq( \frac{\lambda_{M}}{m} +2\epsilon_{0}) | s |^{m\gamma}-C_0,\;\forall (x,s) \in \O \times\mathbb{R},
\end{eqnarray}
and
\begin{eqnarray}\label{e44'}
\lambda_{M}\int_\O|\varphi|^{m\gamma} \leq  \widehat{M}(\|\varphi\|^m)\leq (\lambda_{M}+m\epsilon_{0})\int_\O|\varphi|^{m\gamma}.
\end{eqnarray}
 Set $v=t\varphi,\;t\geq 1$ and using \eqref{e44}, we obtain
\begin{eqnarray}\label{e33}
I(v)&\leq& \frac{1}{m}\widehat{M}(t^m\|\varphi\|^m)
 -( \frac{\lambda_{M}}{m} +2\epsilon_{0})t^{m\gamma} \int_\Omega |\varphi |^{m\gamma}+C_0 |\Omega|\nonumber\\
&\leq& \frac{1}{m}\left(\dfrac{\widehat{M}(t^m\|\varphi\|^m)}{t^{m\gamma}}
-(\lambda_{M}+m\epsilon_{0})
\displaystyle\int_\O|\varphi|^{m\gamma}\right)t^{m\gamma} -\epsilon_0t^{m\gamma}\displaystyle\int_\O|\varphi|^{m\gamma}+C_0 |\Omega|.
\end{eqnarray}
     Using now \eqref{MMM} with $\tau_1=\|\varphi\|^m$, we obtain \begin{eqnarray}\label{ttm}
  \dfrac{\widehat{M}(t^m\|\varphi\|^m)}{ t^{m\gamma}}\leq \widehat{M}(\|\varphi\|^m),\; \forall t\geq 1.\end{eqnarray}
   So, from \eqref{e44'} and \eqref{e33}, we derive
  \begin{eqnarray*}\label{e333}
I(t\varphi) \leq \frac{1}{m}\left( \widehat{M}(\|\varphi\|^m) -(\lambda_{M}+m\epsilon_{0})\displaystyle\int_\O|\varphi|^{m\gamma}\right)t^{m\gamma} -\epsilon_0t^{m\gamma}\displaystyle\int_\O|\varphi|^{m\gamma}+C_0 |\Omega| \leq   -\epsilon_0t^{m\gamma}\displaystyle\int_\O|\varphi|^{m\gamma}+C_0 |\Omega|.\nonumber
\end{eqnarray*}
Choose $t$ large enough, we deduce that $I(v)<0$. In conclusion, $I$ satisfies the mountain pass geometry which ends the proof of Theorem \ref{th24}. \qed

\section*{Acknowledgements}
The authors wish to thank Professor Dong Ye for stimulating discussions on the subject. Also, the second author would like to express his deepest gratitude to the International Centre for Theoretical Physics (ICTP), Trieste, Italie for providing him with an excellent atmosphere for doing this work.

\end{document}